\documentclass[a4paper,final]{article}

\usepackage[T1]{fontenc}
\usepackage[utf8]{inputenc}

\usepackage{amsmath,amssymb,amsthm}

\usepackage{mathtools}

\newtheorem{teo}{Theorem}

\newtheorem{lem}[teo]{Lemma}
\newtheorem{cor}[teo]{Corollary}

\theoremstyle{definition}
\newtheorem{defi}[teo]{Definition}
\newtheorem{rem}[teo]{Remark}
\newenvironment{dimo}
{\begin{proof}}
{\end{proof}}

\usepackage{graphicx}

\newcommand{\matRP}{\mathbb{RP}}
\newcommand{\matN}{\mathbb{N}}
\newcommand{\matR}{\mathbb{R}}
\newcommand{\matS}{\mathbb{S}}
\newcommand{\matB}{\mathbb{B}}
\newcommand{\matT}{\mathbb{T}}
\newcommand{\matK}{\mathbb{K}}
\newcommand{\matA}{\mathbb{A}}
\newcommand{\matM}{\mathbb{M}}
\newcommand{\matL}{\mathbb{L}}

\newcommand{\calC}{{\mathcal C}}
\newcommand{\calR}{{\mathcal R}}
\newcommand{\calT}{{\mathcal T}}

\newcommand{\ptwoirred}{$\mathbb{P}^2$-irreducible}

\newcommand{\Lthreeone}{\matL_{3,1}}
\newcommand{\Lfourone}{\matL_{4,1}}
\newcommand{\twoRPtwo}{2\times\matRP^2}

\newcommand{\timtil}{\begin{picture}(12,12)
\put(2,0){$\times$}\put(2,4.5){$\sim$}\end{picture}}

\newcommand{\co}{\colon\thinspace}

\begin{document}

\title{A complexity of compact 3-manifold via immersed surfaces}

\author{\textsc{Gennaro Amendola}}

\maketitle

\begin{abstract}
We define an invariant, which we call surface-complexity, of compact 3-manifolds by means of Dehn surfaces.
The surface-complexity is a natural number measuring how much the manifold is complicated.
We prove that it fulfils interesting properties: it is subadditive under connected sum and finite-to-one on \ptwoirred\ and boundary-irreducible manifolds without essential annuli and M\"obius strips.
Moreover, for these manifolds, it equals the minimal number of cubes in a cubulation of the manifold, except for the sphere, the ball, the projective space and the lens space $\Lfourone$, which have surface-complexity zero.
We will also give estimations of the surface-complexity by means of ideal triangulations and Matveev complexity.
\end{abstract}

3-manifold, complexity, immersed surface, (ideal) cubulation.


\section*{Introduction}

The problem of filtering compact 3-manifolds in order to study them systematically has been approached by many mathematicians.
The aim is to find a function from the set of compact 3-manifolds to the set of numbers.
The number associated to a 3-manifold should be a measure of how much the manifold is complicated.
For instance, for compact surfaces, this can be achieved by means of genus.
For compact 3-manifolds, many possible functions have been found: e.g.~the Heegaard genus, the Gromov norm, the Matveev complexity.

Each of these functions have features that can be used to study (classes of) compact 3-manifolds.
For instance, they are additive under connected sum.
However, some of them have drawbacks.
The Heegaard genus and the Gromov norm are not finite-to-one, while the Matveev complexity is.
Hence, in order to carry out a classification process, the latter one is more suitable than the former ones.
The Matveev complexity is also a natural measure of how much the manifold is complicated, because if a 3-manifold is \ptwoirred, boundary-irreducible, without essential annuli and M\"obius strips, and different from the ball $\matB^3$, the sphere $\matS^3$, the projective space $\matRP^3$ and the Lens space $\Lthreeone$, then its Matveev complexity is the minimal number of tetrahedra in an ideal triangulation of its (the Matveev complexity of $\matB^3$, $\matS^3$, $\matRP^3$ and $\Lthreeone$ is $0$).
Such functions could also be tools to give proofs by induction.
For instance, the Heegaard genus was used by Rourke to prove by induction that every closed orientable 3-manifold is the boundary of a compact orientable 4-manifold~\cite{Rourke}.

For closed 3-manifolds, one of this function is the surface-complexity, defined in~\cite{Amendola:surf_compl} by means of triple points of images of transverse immersions of closed surfaces that divide the manifold into balls.
It is a function from the set of closed 3-manifolds to the set of natural numbers, which fulfils some properties: it is subadditive under connected sum, and, in the \ptwoirred\ case, it is finite-to-one and it equals the minimal number of cubes in a cubulation of the manifold.
Analogous interesting definitions (which inspired the surface-complexity one) are the Montesinos complexity and the triple point spectrum, given by Vigara~\cite{Vigara:calculus} and studied by Lozano and Vigara~\cite{Lozano-Vigara:subadditivity,Lozano-Vigara:representing,Lozano-Vigara:triple_point_spectrum}.
The definition of surface-complexity is similar to that of Montesinos complexity, but it has the advantage of being more flexible, allowing to prove the properties listed above.

The aim of this paper is to generalise the definition of surface-complexity to the compact case, to generalise the properties holding in the closed case and to give bounds.
We plan to give a complete list of compact 3-manifolds with surface-complexity one in a subsequent paper, as done for the orientable closed case in~\cite{Amendola:sc1}, see also~\cite{Korablev-Kazakov:cubic_complexity_two}.

We now sketch out the definition and the results of this paper.
The surface-complexity of a connected and compact 3-manifold will be defined by means of the notion of {\em quasi-filling Dehn surface}, i.e.~the image of a transverse immersion of a closed surface to which the manifold, possibly after removing some open disjoint balls, collapses.
This generalises the notion of quasi-filling Dehn surface given for closed 3-manifolds in~\cite{Amendola:surf_compl}, which in the end turns out to be a weaker version of the notion of filling Dehn surface given by Montesinos~\cite{Montesinos}.
It is worth noting that (quasi-)filling Dehn surfaces are not yet studied as the other structures (e.g.~Heegaard splittings, triangulations, spines) used to define the other functions described above (there are only partial results~\cite{Vigara:calculus,Amendola:surf_inv,Vigara:lifting,Lozano-Vigara:representing}).

\begin{itemize}
\item[{\bf Definition}] The {\em surface-complexity} $sc(M)$ of a connected and compact 3-manifold $M$ is the minimal number of triple points of a quasi-filling Dehn surface of $M$.
\end{itemize}
For the closed case, this notion differs from Montesinos complexity in two points: firstly quasi-fillingness is used instead of fillingness, secondly any closed surface is allowed as the domain of the transverse immersion instead of the sphere only.
Therefore the surface-complexity could be also called {\em weak-Montesinos complexity}, and for closed 3-manifolds the inequality
$sc(M)\leqslant mc(M)$
obviously holds, where $mc(M)$ denotes the Montesinos complexity of $M$.
In this paper, we have decided to stick to ``surface-complexity'' to stress the similarities with the closed case studied in~\cite{Amendola:surf_compl}.

We will prove the following three properties.
\begin{itemize}
\item[{\bf Finiteness}]
For any integer $c$ there exists only a finite number of connected, compact, \ptwoirred\ and boundary-irreducible 3-manifolds without essential annuli and M\"obius strips that have surface-complexity $c$.
\item[{\bf Naturalness}]
The surface-complexity of a connected, compact, \ptwoirred\ and boundary-irreducible 3-manifold without essential annuli and M\"obius strips, different from $\matS^3$, $\matB^3$, $\matRP^3$ and $\Lfourone$, is equal to the minimal number of cubes in an ideal cubulation of $M$.
The surface-complexity of $\matS^3$, $\matB^3$, $\matRP^3$ and $\Lfourone$ is zero.
\item[{\bf Subadditivity}]
The surface-complexity of the connected sum and of the boundary connected sum of connected and compact 3-manifolds is less than or equal to the sum of their surface-complexities.
\end{itemize}
The naturalness property will follow from the features of {\em minimal} quasi-filling Dehn surfaces of connected, compact, \ptwoirred\ and boundary-irreducible 3-manifolds without essential annuli and M\"obius strips, where minimal means ``with a minimal number of triple points''.
A quasi-filling Dehn surface of a connected and compact 3-manifold is called {\em filling} if it is cell-decomposed by its singularities.
The cell-decomposition dual to a filling Dehn-surface is actually an ideal cubulation of the manifold.
Hence, in order to prove the naturalness property, we will prove that every connected, compact, \ptwoirred\ and boundary-irreducible 3-manifold without essential annuli and M\"obius strips, different from $\matS^3$, $\matB^3$, $\matRP^3$ and $\Lfourone$, has a minimal filling Dehn surface.
We point out that not all minimal quasi-filling Dehn surfaces of connected, compact, \ptwoirred\ and boundary-irreducible 3-manifolds without essential annuli and M\"obius strips are indeed filling.
However, they can be all constructed starting from filling ones (except for $\matS^3$, $\matB^3$, $\matRP^3$ and $\Lfourone$, for which non-filling ones must be used) and applying a simple move, we will call {\em bubble-move}.

The surface-complexity is related to the Matveev complexity.
Indeed, if $M$ is a connected, compact, \ptwoirred\ and boundary-irreducible 3-manifold without essential annuli and M\"obius strips, different from $\Lthreeone$ and $\Lfourone$, the double inequality
$\frac{1}{8}c(M) \leqslant sc(M) \leqslant 4c(M)$
holds, where $c(M)$ denotes the Matveev complexity of $M$.

The two inequalities above give also estimates of the surface-complexity.
In general, an exact calculation of the surface-complexity is very difficult, however it is relatively easy to estimate it.
More precisely, it is quite easy to give upper bounds for it, because constructing a quasi-filling Dehn surface of the manifold with the appropriate number of triple points suffices.
With this technique, we will give an upper bound for the surface-complexity of a connected and compact 3-manifold starting from an ideal triangulation of its.
The problem of proving the sharpness of this bound is usually difficult (see for instance \cite{Lozano-Vigara:triple_point_spectrum} where the triple point spectrum of $\matS^3$ and $\matS^2\times\matS^1$ has been computed very subtly).

In the Appendix we will give a brief description of what happens in the 2-dimensional case.
We plan to cope with the 4-dimensional case in a subsequent paper.

\paragraph{Acknowledgements}
We would like to thank Bruno Martelli for the useful discussions on the Matveev complexity.

\section{Definitions}

Throughout this paper, all 3-manifolds are assumed to be connected and compact.
By $M$, we will always denote such a (connected and compact) 3-manifold.
Using the {\em Hauptvermutung}, we will freely intermingle the differentiable,
piecewise linear and topological viewpoints.

The manifold $M$ is called {\em \ptwoirred} if every sphere in $M$ bounds a ball and $M$ does not contain any two-sided embedded projective plane.
The manifold $M$ is called {\em boundary-irreducible} if for every proper disc $D\subset M$ the curve $\partial D$ bounds a disc in $\partial M$.
Suppose $M$ is \ptwoirred\ and boundary-irreducible, then a proper annulus $A\subset M$ is called {\em essential} if $A$ does not cut off from $M$ either a cylinder (a ball intersecting $\partial M$ in two discs) or a solid torus (intersecting $\partial M$ in another annulus), and a proper M\"obius strip $A\subset M$ is called {\em essential} if $A$ does not cut off from $M$ a solid Klein bottle (intersecting $\partial M$ in another M\"obius strip).
The definitions of essential annulus and M\"obius strip are more general, but in a \ptwoirred\ and boundary-irreducible $M$ these definitions are equivalent to the general ones, and for the purpose of this paper we do not need the general ones.

In what follows, we will need to remove some open balls from the manifold $M$ such that their closures and the components of the boundary of $M$ are pairwise disjoint.
In such a case, the manifold obtained will be denoted by $\dot M$.
Note that the boundary of $\dot M$ is made up of $\partial M$ and some spheres.
We allow the case where some components of $\partial M$ are spheres, and in this case they are regarded to be different with respect to the spheres in $\partial\dot M\setminus\partial M$.
If instead $M$ is closed, the boundary of $\dot M$ is made up of the boundary of the removed balls only.
With a slight abuse of notation, we will call $\dot M$ to be {\em a punctured $M$}.

\paragraph{Dehn surfaces}
A subset $\Sigma$ of $M$ is said to be a {\em Dehn surface of $M$}~\cite{Papa}
if there exists an abstract (possibly non-connected) closed surface $S$ and a transverse immersion $f\co S\to M$ such that $\Sigma = f(S)$.

\begin{figure}
  \centerline{
  \begin{tabular}{ccc}
    \begin{minipage}[c]{3.5cm}{\small{\begin{center}
        \includegraphics{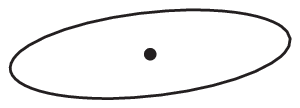}
      \end{center}}}\end{minipage} &
    \begin{minipage}[c]{3.5cm}{\small{\begin{center}
        \includegraphics{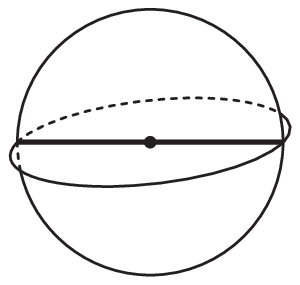}
      \end{center}}}\end{minipage} &
    \begin{minipage}[c]{3.5cm}{\small{\begin{center}
        \includegraphics{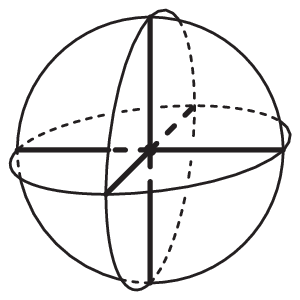}
      \end{center}}}\end{minipage} \\
    \begin{minipage}[t]{3.5cm}{\small{\begin{center}
        Simple\\point
      \end{center}}}\end{minipage} &
    \begin{minipage}[t]{3.5cm}{\small{\begin{center}
        Double\\point
      \end{center}}}\end{minipage} &
    \begin{minipage}[t]{3.5cm}{\small{\begin{center}
        Triple\\point
      \end{center}}}\end{minipage}
  \end{tabular}}
  \caption{Neighbourhoods of points (marked by thick dots) of a Dehn surface.}
  \label{fig:neigh_Dehn_surf}
\end{figure}

Let us fix for a while $f\co S\to M$ a transverse immersion (hence,
$\Sigma = f(S)$ is a Dehn surface of $M$). By transversality, the number
of pre-images of a point of $\Sigma$ is 1, 2 or 3; so there are three types of
points in $\Sigma$, depending on this number; they are called
{\em simple}, {\em double} or {\em triple}, respectively.
Note that the definition of the type of a point does not depend on the particular
transverse immersion $f\co S\to M$ we have chosen. In fact, the
type of a point can be also defined by looking at a regular neighbourhood (in
$M$) of the point, as shown in Fig.~\ref{fig:neigh_Dehn_surf}.
The set of triple points is denoted by $T(\Sigma)$; non-simple points are called {\em singular} and their set is denoted by $S(\Sigma)$.
From now on, in all figures, triple points are always marked by thick dots and the singular set is also drawn thick.

\paragraph{(Quasi-)filling Dehn surfaces}
A Dehn surface $\Sigma$ of $M$ will be called {\em quasi-filling} if a punctured $M$, $\dot M$, collapses to $\Sigma$.
Moreover, a quasi-filling $\Sigma$ is called {\em filling} if its singularities induce a cell-decomposition of $\Sigma$; more precisely,
\begin{itemize}
\item $T(\Sigma) \neq \emptyset$,
\item $S(\Sigma) \setminus T(\Sigma)$ is made up of intervals (called {\em
    edges}),
\item $\Sigma \setminus S(\Sigma)$ is made up of discs (called {\em regions}).
\end{itemize}

Since $M$ is connected, any quasi-filling Dehn surface $\Sigma$ of $M$ is connected.
By construction, the complement of $\Sigma$ in $M$ is made up of components of two types:
\begin{itemize}
\item
the components intersecting $\partial M$, which are all together isomorphic to $\partial M\times[0,1)$, with $\partial M$ corresponding to $\partial M\times\{0\}$;
\item
the components intersecting $\partial\dot M\setminus\partial M$, which are open balls.
\end{itemize}
Moreover, note that $\Sigma$ is also a quasi-filling Dehn surface of $\dot M$, with all components of $\dot M\setminus\Sigma$ being of the former type.
Finally, note that any small regular neighbourhood $\calR(\Sigma)$ of $\Sigma$ in $M$ is isomorphic to $\dot M$.

\begin{rem}
If $M$ is closed, a Dehn surface $\Sigma$ of $M$ is quasi-filling if $M \setminus \Sigma$ is made up of balls.
Therefore, this notion generalises the notion of (quasi-)filling Dehn surfaces of closed orientable 3-manifolds to compact ones~\cite{Montesinos,Amendola:surf_inv}.
\end{rem}

\begin{rem}\label{rem:Sigma_surf}
Suppose a quasi-filling $\Sigma$ of $M$ is a surface (i.e.~$S(\Sigma)=\emptyset$).
Then, $\calR(\Sigma)\cong\dot M$ is an $I$-bundle over $\Sigma$.
Since the boundary of $\dot M$ is the union of $\partial M$ and spheres, and since $\Sigma$ is connected, we have four cases as described in Table~\ref{tab:surface}.
\begin{table}
	\centerline{\begin{tabular}{cccc}
	\phantom{\Big|}$\Sigma$ & $M$ & $\dot M$ & $\#$ of removed balls \\ \hline
	\phantom{\Big|}$\matS^2$ & $\matS^3$ & $\matS^2\times[-1,1]$ & $2$ \\
	\phantom{\Big|}$\matS^2$ & $\matB^3$ & $\matS^2\times[-1,1]$ & $1$ \\
	\phantom{\Big|}$\matRP^2$ & $\matRP^3$ & $\matRP^2\timtil[-1,1]$ & $1$ \\
	\phantom{\Big|}any & $I$-bundle over $\Sigma$  & $I$-bundle over $\Sigma$ & $0$
	\end{tabular}}
	\caption{The cases where $\Sigma$ is a surface.}
	\label{tab:surface}
\end{table}
\end{rem}

Let us give some other examples.
Two projective planes intersecting along a loop non-trivial in both of them, which will be called {\em double projective plane} and denoted by $\twoRPtwo$, form a quasi-filling Dehn surface (without triple points) of $\matRP^3$ or $\matRP^3$ with up to 2 balls removed.
A sphere intersecting a torus (resp.~a Klein bottle) along a loop is a quasi-filling Dehn surface (without triple points) of $\matS^2\times\matS^1$ (resp.~$\matS^2\timtil\matS^1$) or $\matS^2\times\matS^1$ (resp.~$\matS^2\timtil\matS^1$) with up to 2 balls removed.
A sphere self-intersecting in a loop is a quasi-filling Dehn surface (without triple points) of a solid torus or of a solid Klein bottle (depending on the orientation of the self-intersection) with up to 2 balls removed.
The {\em quadruple hat} (i.e.~a disc whose boundary is glued four times along a circle) is a quasi-filling Dehn-surface (without triple points) of the lens-space $\Lfourone$ or of $\Lfourone$ with 1 ball removed.
If we identify the sphere $\matS^3$ with $\matR^3\cup\{\infty\}$, the three coordinate planes in $\matR^3$, with $\infty$ added, form a filling Dehn surface (with two triple points: $(0,0,0)$ and $\infty$) of $\matS^3$ or $\matS^3$ with up to 8 balls removed.

\paragraph{Ideal cubulations and duality}

Consider the topological space $\widehat M$ obtained from $\dot M$ by collapsing each boundary component (of $\dot M$) to a point.
Note that the complement of the points corresponding to the boundary components of $\dot M$ can be identified with the interior of $\dot M$, and that $\widehat M$ can be obtained also from $M$ by collapsing each boundary component (of $M$) to a point.

\begin{figure}
  \centerline{\includegraphics{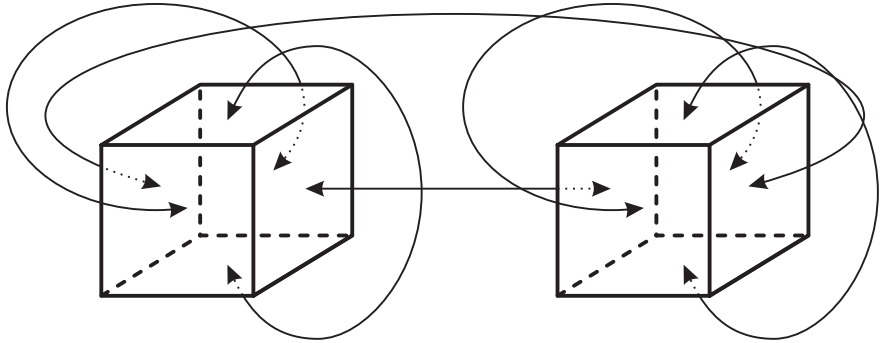}}
  \caption{A cubulation of the 3-dimensional torus $\matS^1\times\matS^1\times\matS^1$ with two cubes (the identification of each pair of faces is the obvious one, i.e.~the one without twists).}
  \label{fig:cubul_example}
\end{figure}

\begin{figure}
  \centerline{\includegraphics{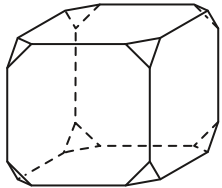}}
  \caption{A cube with small tetrahedra removed near the vertices.}
  \label{fig:cube_tetra_removed}
\end{figure}

An {\em ideal cubulation} $\calC$ of $M$ is a cell-decomposition of $\widehat M$ such that
\begin{itemize}
\item the set of 0-cells (called {\em vertices}) is the set of the points corresponding to the boundary components of $\dot M$,
\item each 2-cell (called a {\em face}) is glued along 4 edges,
\item each 3-cell (called a {\em cube}) is glued along 6 faces arranged like the boundary of a cube.
\end{itemize}
Note that self-adjacencies and multiple adjacencies are allowed.
In Fig.~\ref{fig:cubul_example} we have shown a cubulation of the 3-dimensional torus $\matS^1\times\matS^1\times\matS^1$ with two cubes (a closed manifold).
There are two types of vertices in $\calC$: those corresponding to the boundary components of $M$ and those corresponding to the spheres of $\partial\dot M\setminus\partial M$.
We will call {\em ideal} the former ones and {\em finite} the latter ones.
Note that the whole $\dot M$ can be obtained from $\calC$ by removing small tetrahedra in each cube near the vertices (see Fig.~\ref{fig:cube_tetra_removed}), that the complement of the ideal vertices can be identified with the interior of $M$, and that $M$ can be obtained from $\calC$ by removing small tetrahedra in each cube near the ideal vertices.
It is worth noting that $\calC$ does not distinguish the ideal vertices corresponding to spheres from the finite ones.
Note also that ours is a slight abuse of notation, indeed in the classical meaning an ideal polyhedron has only ideal vertices.

\begin{figure}
  \centerline{\includegraphics{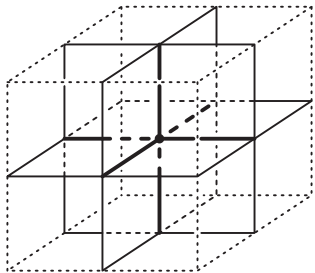}}
  \caption{Local behaviour of duality.}
  \label{fig:cube_to_surf}
\end{figure}

The following construction is well-known for cubulations (see~\cite{Aitchison-Matsumotoi-Rubinstein, Funar, Babson-Chan}, for instance), but applies also to ideal cubulations.
Let $\calC$ be an ideal cubulation of $M$.
Consider, for each cube of $\calC$, the three squares shown in Fig.~\ref{fig:cube_to_surf}.
We can suppose that the squares fit together through the faces, so the subset of $M$ obtained by gluing together these squares is a filling Dehn surface $\Sigma$ of $M$.
Conversely, an ideal cubulation $\calC$ of $M$ can be constructed from a filling Dehn surface $\Sigma$ of $M$ by considering an abstract cube for each triple point of $\Sigma$ and by gluing the cubes together along the faces, with the identification of each pair of faces chosen by following the four germs of regions adjacent to the respective edge of $\Sigma$.
The ideal cubulation and the filling Dehn surface constructed in such a way are said to be {\em dual} to each other.

\paragraph{Existence of filling Dehn surfaces and reconstruction of the manifold}
It is by now well-known that every closed $M$ has a filling Dehn-surface (see, for instance,~\cite{Montesinos,Vigara:present,Amendola:surf_inv} and also~\cite{Fenn-Rourke} for the quasi-filling case).
We will adapt the proof of~\cite{Amendola:surf_inv} to the non-closed case.
In what follows, as we have done above for cubulations, with a slight abuse of notation, we use a non classical notion of {\em ideal triangulation}, i.e.~a triangulation $\calT$ of $\widehat M$ whose vertices are the points corresponding to the boundary components of $\dot M$.
We have {\em ideal} and {\em finite} vertices: the former being those corresponding to the boundary components of $M$, the latter being those corresponding to the spheres of $\partial\dot M\setminus\partial M$.

\begin{teo}\label{teo:filling_surf_present}
\begin{itemize}
\item
Each (connected and compact) 3-manifold has a filling Dehn surface.
\item
If $\Sigma_1$ and $\Sigma_2$ are homeomorphic filling Dehn surfaces of (connected and compact) 3-manifolds $M_1$ and $M_2$, respectively, with the same number of spherical boundary components, then $M_1$ and $M_2$ are also homeomorphic.
\end{itemize}
\end{teo}

\begin{dimo}
We start by proving the first point.
Let $\calT$ be an ideal triangulation of a 3-manifold $M$ (all 3-manifolds have ideal triangulations, as shown in~\cite{Matveev:book}, for instance).
Consider, for each tetrahedron of $\calT$, 4 triangles close and parallel to the 4 faces, as shown in Fig.~\ref{fig:tria_to_surf}.
\begin{figure}
  \centerline{\includegraphics{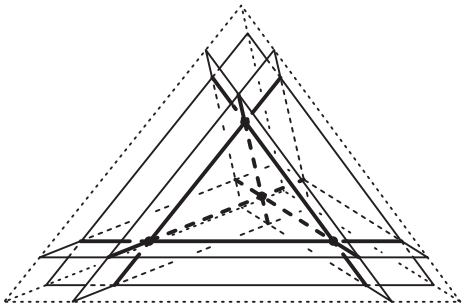}}
  \caption{Construction of a filling Dehn sphere from an ideal triangulation.}
  \label{fig:tria_to_surf}
\end{figure}
We can suppose that the triangles fit together through the faces, so the subset of $M$ obtained by gluing together these triangles is a Dehn surface $\Sigma$ of $M$.
It is very easy to prove that $\Sigma$ is filling, so we leave it to the reader.

We do not give a complete proof of the second point because it is essentially the same as that of Casler for standard spines~\cite{Casler}.
The idea of the proof is the following.
Let $\calC_i$ be the cubulation of $M_i$ dual to $\Sigma_i$, for $i=1,2$.
The cubulations $\calC_i$ are defined unambiguously, because the cubes dual to the triple points and the face identifications of them are defined unambiguously (up to homeomorphism) from the Dehn surfaces $\Sigma_i$.
Since $\Sigma_1$ and $\Sigma_2$ are homeomorphic, the cubulations $\calC_1$ and $\calC_2$ turn out to be isomorphic and hence $\dot{M_1}$ and $\dot {M_2}$ are homeomorphic.
The fact that $M_1$ and $M_2$ have the same number of spherical boundary components easily implies that they are obtained from $\dot{M_1}$ and $\dot {M_2}$, respectively, by filling up the same number of boundary spherical components, so they are homeomorphic.
\end{dimo}

With respect to the construction described in the first part of the proof above, it is worth noting that it is the dual counterpart of the well-known construction consisting in dividing a tetrahedron into 4 cubes~\cite{Shtanko-Shtogrin, Dolbilin-Shtanko-Shtogrin, Funar}.
Note also that the surface $S$ such that $\Sigma = f(S)$ is $\partial\dot M$, because $\Sigma$ can be obtained by starting with a copy of $\partial\dot M$ parallel to $\partial\dot M$ in the interior of $\dot M$ and then moving it away from the boundary.
Finally, note that $\Sigma$ has 4 triple points for each tetrahedron.

\begin{rem}\label{rem:mfld_uniquely_determined}
With respect to the construction described in the second part of the proof above, if we require that the $M_i$'s are \ptwoirred\ and different from the ball $\matB^3$, the $M_i$'s have no boundary spherical component, so we must fill up all spherical boundary components of the $\dot M_i$'s to recover the $M_i$'s, and hence the $M_i$'s are homeomorphic.
\end{rem}

\paragraph{Abstract filling Dehn surfaces}
A filling Dehn surface $\Sigma$ of $M$ is contained in $M$.
However, one can think of it as an abstract cell complex.
For the sake of completeness, we mention that the abstract cell complex $\Sigma$ determines $\dot M$ and $M$ if the number of spherical boundary components is known (and the abstract surface $S$ such that $\Sigma=f(S)$ where $f \co S \rightarrow M$) up to homeomorphism (the proof is the same as that of the second part of Theorem~\ref{teo:filling_surf_present}).
This is not true for quasi-filling Dehn surfaces, e.g.~a surface $\Sigma$ is a quasi-filling Dehn surface of any $I$-bundle over $\Sigma$.

\paragraph{Surface-complexity}
The surface-complexity of $M$ can now be defined as the minimal number of triple points of a quasi-filling Dehn surface of $M$. More precisely, we give the following.
\begin{defi}
The {\em surface-complexity} $sc(M)$ of $M$ is equal to $c$ if $M$ possesses a quasi-filling Dehn surface with $c$ triple points and has no quasi-filling Dehn surface with less than $c$ triple points.
In other words, $sc(M)$ is the minimum of $|T(\Sigma)|$ over all quasi-filling Dehn surfaces $\Sigma$ of $M$.
\end{defi}

We will classify the 3-manifolds having surface-complexity zero in the following section.
At the moment, we can only say that $\matS^3$, $\matB^3$, $\matRP^3$, $\matS^2\times\matS^1$, $\matS^2\timtil\matS^1$, $\Lfourone$, the $I$-bundles and some of these manifolds with a certain number of balls removed have surface-complexity 0, because we have seen above that they have quasi-filling Dehn surfaces without triple points.

\section{Minimality and finiteness}

A quasi-filling Dehn surface $\Sigma$ of $M$ is called {\em minimal} if it has a minimal number of triple points among all quasi-filling Dehn surfaces of $M$, i.e.~$|T(\Sigma)|=sc(M)$.

\begin{teo}\label{teo:minimal_filling}
Let $M$ be a (connected and compact) \ptwoirred\ and boundary-irreducible 3-manifold without essential annuli and M\"obius strips.
\begin{itemize}
\item
If $sc(M)=0$, the manifold $M$ is the sphere $\matS^3$, the ball $\matB^3$, the projective space $\matRP^3$ or the lens space $\Lfourone$.
\item
If $sc(M)>0$, the manifold $M$ has a minimal filling Dehn surface.
\end{itemize}
\end{teo}

\begin{dimo}
Let $\Sigma$ be a minimal quasi-filling Dehn surface of $M$.
If we have $S(\Sigma)=\emptyset$ (i.e.~$\Sigma$ is a surface), by virtue of Remark~\ref{rem:Sigma_surf}, we have that $M$ is $\matS^3$, $\matB^3$, $\matRP^3$ or an $I$-bundle.
The last case cannot occur, because an $I$-bundle is not \ptwoirred, or contains essential annuli or M\"obius strips.

Then, we suppose $S(\Sigma)\neq\emptyset$.
We will first prove that $M$ has a quasi-filling Dehn surface $\Sigma'$ such that $\Sigma' \setminus S(\Sigma')$ is made up of discs or $\Sigma'$ is a surface.
Suppose there exists a component $C$ of $\Sigma \setminus S(\Sigma)$ that is not a disc.
$C$ contains a non-trivial orientation preserving (in $C$) simple closed curve $\gamma$.
Consider a strip $F$ properly embedded in $\dot M$ such that $F\cap\Sigma=\gamma$.
The strip $F$ is either an annulus or a M\"obius strip depending on whether $\gamma$ is orientation preserving in $M$ or not.
Since $M\setminus\dot M$ is made up of balls, we can fill up $\partial F$ with two, one or zero discs in $M\setminus\dot M$, depending on whether the number of components of $F\cap(\partial\dot M\setminus\partial M)$ are two, one or zero, respectively.
We then get a surface $\overline{F}$ properly embedded in $M$.
If $F$ is an annulus, $\overline{F}$ is a sphere, a disc or $F$ itself.
If $F$ is an M\"obius strip, $\overline{F}$ is a projective plane or $F$ itself.
Therefore, we have five cases to analyse.
\begin{itemize}
\item
If $\overline{F}$ is a sphere, it bounds a ball (in $M$), say $B$, because $M$ is \ptwoirred.
Since $\Sigma\cap\partial B=\gamma$ is a simple closed curve, we can replace the portion of $\Sigma$ contained in $B$ with a disc, getting a new quasi-filling Dehn surface of $M$.
\item
If $\overline{F}$ is a disc, it cuts off a ball from $M$ (intersecting $\partial M$ in another disc), say $B$, because $M$ is \ptwoirred\ and boundary-irreducible.
Since $\Sigma\cap\partial B=\gamma$ is a simple closed curve, we can replace the portion of $\Sigma$ contained in $B$ with a disc, getting a new quasi-filling Dehn surface of $M$.
\item
If $\overline{F}$ is the annulus $F$, it cuts off from $M$ either a cylinder (intersecting $\partial M$ in two discs), say $B$, or a solid torus (intersecting $\partial M$ in another annulus), say $T$, because $M$ is \ptwoirred, boundary-irreducible and without essential annuli.
In the former case, since $\Sigma\cap\partial B=\gamma$ is a simple closed curve, we can replace the portion of $\Sigma$ contained in $B$ with a disc, getting a new quasi-filling Dehn surface of $M$.
The latter case instead cannot occur, because the core of $\overline{F}$ would be a longitude of the solid torus $T$ and would be null-homologous in $T$ (indeed $\Sigma\cap T$ is a 2-cycle whose boundary is the core of $\overline{F}$), a contradiction.
\item
If $\overline{F}$ is a projective plane, it is two-sided (because $F$ is), which is not possible because $M$ is \ptwoirred.
\item
If $\overline{F}$ is the M\"obius strip $F$, it cuts off from $M$ a solid Klein bottle (intersecting $\partial M$ in another M\"obius strip), say $K$, because $M$ is \ptwoirred, boundary-irreducible and without essential M\"obius strips.
This case cannot occur, because the core of $\overline{F}$ would be a longitude of the Klein bottle $K$ and would be null-homologous in $K$ (indeed $\Sigma\cap K$ is a 2-cycle whose boundary is the core of $\overline{F}$), a contradiction.
\end{itemize}
In all admissible cases the Euler characteristic of the component of $\Sigma \setminus S(\Sigma)$ containing $\gamma$ has increased, no new non-disc component has been created and the number of triple points is not increased.
Hence, by repeatedly applying this procedure, we eventually get a quasi-filling Dehn surface, say $\Sigma'$, of $M$ such that $\Sigma' \setminus S(\Sigma')$ is made up of discs or $\Sigma'$ is a surface.

If $\Sigma'$ is a surface, analogously to what said at the beginning of the proof for $\Sigma$, by virtue of Remark~\ref{rem:Sigma_surf}, we have that $M$ is $\matS^3$, $\matB^3$ or $\matRP^3$.

Therefore, we eventually consider the case where $\Sigma' \setminus S(\Sigma')$ is made up of discs.
Since $\Sigma'$ is connected, also $S(\Sigma')$ is connected.
If we have $sc(M)>0$ (i.e.~$T(\Sigma')$ is not empty), $S(\Sigma') \setminus T(\Sigma')$ cannot contain circles and hence $\Sigma'$ is filling (i.e.~$M$ has a minimal filling Dehn surface).
Otherwise, if we have $sc(M)=0$ (i.e.~$T(\Sigma')$ is empty), $S(\Sigma')$ is made up of one circle.
Since $\Sigma' \setminus S(\Sigma')$ is made up of discs, the Dehn surface $\Sigma'$ is completely determined by the regular neighbourhood of $S(\Sigma')$ in $\Sigma'$.
This neighbourhood depends on how the germs of disc are interchanged along the curve $S(\Sigma')$.
Since $M$ is \ptwoirred, only three cases must be taken into account for $\Sigma'$ (up to symmetry):
\begin{itemize}
\item two spheres intersecting along the circle $S(\Sigma')$, which form a Dehn surface of $\matS^3$ or $\matB^3$;
\item the double projective plane $\twoRPtwo$, which is a Dehn surface of $\matRP^3$;
\item the four-hat, which is a Dehn surface of $\Lfourone$.
\end{itemize}
This concludes the proof.
\end{dimo}

Since there is a finite number of filling Dehn surfaces having a fixed number of triple points and each of them is the Dehn surface of one \ptwoirred\ manifold (two for the $\matS^3$/$\matB^3$ case -- see Remark~\ref{rem:mfld_uniquely_determined}), and since there are only a finite number of \ptwoirred\ and boundary-irreducible 3-manifolds without essential annuli and M\"obius strips with surface-complexity 0, we have the following result.

\begin{cor}
For any integer $c$ there exists only a finite number of (connected and compact) \ptwoirred\ and boundary-irreducible 3-manifolds without essential annuli and M\"obius strips that have surface-complexity $c$.
\end{cor}

By means of the duality between filling Dehn surfaces and ideal cubulations we have the following result.

\begin{cor}
The surface-complexity of a (connected and compact) \ptwoirred\ and boundary-irreducible 3-manifold without essential annuli and M\"obius strips, different from $\matS^3$, $\matB^3$, $\matRP^3$ and $\Lfourone$, is equal to the minimal number of cubes in an ideal cubulation of $M$.
\end{cor}

\subsection{Minimal quasi-filling Dehn surfaces}

\begin{figure}
  \centerline{\includegraphics{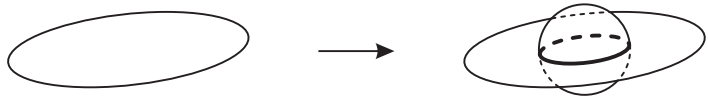}}
  \caption{Bubble-move.}
  \label{fig:bubble_move}
\end{figure}

Theorem~\ref{teo:minimal_filling} states that, under some hypotheses, $M$ has a minimal filling Dehn surface, but not all minimal quasi-filling Dehn surfaces of such an $M$ are indeed filling.
However, we will see now that they can be all constructed starting from filling ones and applying a simple move.
The move acts on quasi-filling Dehn surfaces near a simple point as shown in Fig.~\ref{fig:bubble_move} and it is called a {\em bubble-move}.
Note that the result of applying a bubble-move to a quasi-filling Dehn surface of $M$ is a quasi-filling Dehn surface of $M$, but the result of applying a bubble-move to a filling Dehn-surface is not filling any more.
Note also that the bubble-move increases (by two) the number of balls in $M\setminus\dot M$.
If a quasi-filling Dehn surface $\Sigma$ is obtained from a quasi-filling Dehn surface $\overline{\Sigma}$ by repeatedly applying bubble-moves (even zero), we will say that $\Sigma$ {\em is derived from} $\overline{\Sigma}$.
Obviously, if $\Sigma$ is a quasi-filling Dehn surface of $M$ and is derived from $\overline{\Sigma}$, also $\overline{\Sigma}$ is a quasi-filling Dehn surface of $M$.
Eventually, it is worth noting (even if not useful for the purpose of this paper) that if we have in $\Sigma$ a configuration as in Fig.~\ref{fig:bubble_move}-right, we can apply an inverse bubble move, removing the sphere from $\Sigma$, only if the ball bounded by the sphere in Fig.~\ref{fig:bubble_move}-right is actually a ball in $M$, i.e.~it contains two balls of $M\setminus\dot M$.

We will need the following result, which has been proved in~\cite[Lemma~2.3]{Amendola:surf_compl}.

\begin{lem}\label{lem:all_minimal_filling_sphere}
Let $\Sigma$ be a minimal quasi-filling Dehn surface of the sphere $\matS^3$ and let $D$ be a closed disc contained in $\Sigma \setminus S(\Sigma)$.
Then $\Sigma$ is derived from a sphere $\matS^2$ by means of bubble-moves not involving $D$.
\end{lem}

Theorem~\ref{teo:minimal_filling} can be improved by means of a slightly subtler analysis.

\begin{teo}
Let $\Sigma$ be a minimal quasi-filling Dehn surface of a (connected and compact) \ptwoirred\ and boundary-irreducible 3-manifold $M$ without essential annuli and M\"obius strips.
\begin{itemize}
\item
If $sc(M)=0$, one of the following holds:
\begin{itemize}
\item
$M$ is the sphere $\matS^3$ or the ball $\matB^3$, and $\Sigma$ is derived from the sphere $\matS^2$;
\item
$M$ is the projective space $\matRP^3$, and $\Sigma$ is derived from the projective plane $\matRP^2$ or from the double projective plane $\twoRPtwo$;
\item
$M$ is the lens space $\Lfourone$, and $\Sigma$ is derived from the four-hat.
\end{itemize}
\item
If $sc(M)>0$, the Dehn surface $\Sigma$ is derived from a minimal filling Dehn surface of $M$.
\end{itemize}
\end{teo}

\begin{dimo}
The scheme of the proof is the same as that of Theorem~\ref{teo:minimal_filling}.
Hence, we will often refer to the proof of Theorem~\ref{teo:minimal_filling} also for notation.

Let $\Sigma$ be a minimal quasi-filling Dehn surface of $M$.
If we have $S(\Sigma)=\emptyset$ (i.e.~$\Sigma$ is a surface), by virtue of Remark~\ref{rem:Sigma_surf}, we have three cases:
\begin{itemize}
\item
$\Sigma$ is $\matS^2$ and $M$ is $\matS^3$ or $\matB^3$;
\item
$\Sigma$ is $\matRP^2$ and $M$ is $\matRP^3$;
\item
$\Sigma$ is a surface and $M$ is an $I$-bundle, which is not \ptwoirred, or contains essential annuli or M\"obius strips, so this case cannot occur.
\end{itemize}

Then, we suppose $S(\Sigma)\neq\emptyset$.
We will first prove that $\Sigma$ is derived from a (minimal) quasi-filling Dehn surface $\Sigma'$ of $M$ such that either $\Sigma' \setminus S(\Sigma')$ is made up of discs or $\Sigma'$ is a surface.
Suppose there exists a component $C$ of $\Sigma \setminus S(\Sigma)$ that is not a disc.
$C$ contains a non-trivial orientation preserving (in $C$) simple closed curve $\gamma$.
As done in the proof of Theorem~\ref{teo:minimal_filling}, we consider a strip $F$ properly embedded in $\dot M$ such that $F\cap\Sigma=\gamma$, and fill up $\partial F$ with two, one or zero discs in $M\setminus\dot M$, getting a surface $\overline{F}$ properly embedded in $M$.
Therefore, we have five cases to analyse.
\begin{itemize}
\item
If $\overline{F}$ is a sphere, it bounds a ball (in $M$), say $B$, because $M$ is \ptwoirred.
Consider $\Sigma_B = \Sigma \cap B$ and $\Sigma_M = \Sigma \setminus B$.
If we fill up $\Sigma_M$ with a disc by gluing it along $\gamma$, we obtain a minimal quasi-filling Dehn surface $\Sigma_M'$ of $M$.
Analogously, if we fill up $\Sigma_B$ with a disc (say $D$) by gluing it along $\gamma$, we obtain a minimal quasi-filling Dehn surface $\Sigma_B'$ of $\matS^3$.
By virtue of Lemma~\ref{lem:all_minimal_filling_sphere}, $\Sigma_B'$ is derived from a sphere $\matS^2$ by means of bubble-moves not involving $D$.
These moves can be applied to $\Sigma_M'$ because they do not involve $D$ and the result is $\Sigma$, so $\Sigma$ is derived from the quasi-filling Dehn surface $\Sigma_M'$ of $M$.
\item
If $\overline{F}$ is a disc, it cuts off a ball from $M$ (intersecting $\partial M$ in another disc), say $B$, because $M$ is \ptwoirred\ and boundary-irreducible.
Consider $\Sigma_B = \Sigma \cap B$ and $\Sigma_M = \Sigma \setminus B$.
If we fill up $\Sigma_M$ with a disc by gluing it along $\gamma$, we obtain a minimal quasi-filling Dehn surface $\Sigma_M'$ of $M$.
Analogously, if we fill up $\Sigma_B$ with a disc (say $D$) by gluing it along $\gamma$, we obtain a minimal quasi-filling Dehn surface $\Sigma_B'$ of $\matS^3$.
By virtue of Lemma~\ref{lem:all_minimal_filling_sphere}, $\Sigma_B'$ is derived from a sphere $\matS^2$ by means of bubble-moves not involving $D$.
These moves can be applied to $\Sigma_M'$ because they do not involve $D$ and the result is $\Sigma$, so $\Sigma$ is derived from the quasi-filling Dehn surface $\Sigma_M'$ of $M$.
\item
If $\overline{F}$ is the annulus $F$, it cuts off from $M$ either a solid torus (intersecting $\partial M$ in another annulus), say $T$, or a cylinder (intersecting $\partial M$ in two discs), say $B$, because $M$ is \ptwoirred, boundary-irreducible and without essential annuli.
The former case cannot occur (see the proof of Theorem~\ref{teo:minimal_filling}).
The latter case instead can occur.
Consider $\Sigma_B = \Sigma \cap B$ and $\Sigma_M = \Sigma \setminus B$.
If we fill up $\Sigma_M$ with a disc by gluing it along $\gamma$, we obtain a minimal quasi-filling Dehn surface $\Sigma_M'$ of $M$.
Analogously, if we fill up $\Sigma_B$ with a disc (say $D$) by gluing it along $\gamma$, we obtain a minimal quasi-filling Dehn surface $\Sigma_B'$ of $\matS^3$.
By virtue of Lemma~\ref{lem:all_minimal_filling_sphere}, $\Sigma_B'$ is derived from a sphere $\matS^2$ by means of bubble-moves not involving $D$.
These moves can be applied to $\Sigma_M'$ because they do not involve $D$ and the result is $\Sigma$, so $\Sigma$ is derived from the quasi-filling Dehn surface $\Sigma_M'$ of $M$.
\item
The last two cases, where $\overline{F}$ is a projective plane or the M\"obius strip $F$, cannot occur because $M$ is \ptwoirred, boundary-irreducible and without essential M\"obius strips (see the proof of Theorem~\ref{teo:minimal_filling}).
\end{itemize}
In all admissible cases the Euler characteristic of the component of $\Sigma_M' \setminus S\left(\Sigma_M'\right)$ containing $\gamma$ is bigger than that of the corresponding component of $\Sigma \setminus S(\Sigma)$, no new non-disc component has been created and the number of triple points has not changed.
Hence, by repeatedly applying this procedure, we eventually get a (minimal) quasi-filling Dehn surface $\Sigma'$ of $M$ from which $\Sigma$ is derived and such that either $\Sigma' \setminus S(\Sigma')$ is made up of discs or $\Sigma'$ is a surface.

If $\Sigma'$ is a surface, analogously to what said at the beginning of the proof for $\Sigma$, by virtue of Remark~\ref{rem:Sigma_surf}, we have two cases:
\begin{itemize}
\item
$\Sigma'$ is $\matS^2$ and $M$ is $\matS^3$ or $\matB^3$;
\item
$\Sigma'$ is $\matRP^2$ and $M$ is $\matRP^3$.
\end{itemize}

Therefore, we eventually consider the case where $\Sigma' \setminus S(\Sigma')$ is made up of discs.
Since $\Sigma'$ is connected, also $S(\Sigma')$ is connected.
If we have $sc(M)>0$ (i.e.~$T(\Sigma')$ is not empty), $S(\Sigma') \setminus T(\Sigma')$ cannot contain circles and hence $\Sigma'$ is filling (i.e.~$\Sigma$ is derived from a minimal filling Dehn surface of $M$).
Otherwise, if we have $sc(M)=0$ (i.e.~$T(\Sigma')$ is empty), $S(\Sigma')$ is made up of one circle.
Since $\Sigma' \setminus S(\Sigma')$ is made up of discs, the Dehn surface $\Sigma'$ is completely determined by the regular neighbourhood of $S(\Sigma')$ in $\Sigma'$.
This neighbourhood depends on how the germs of disc are interchanged along the curve $S(\Sigma')$.
Since $M$ is \ptwoirred, only three cases must be taken into account for $\Sigma'$ (up to symmetry):
\begin{itemize}
\item two spheres intersecting along the circle $S(\Sigma')$, which form a Dehn surface of $\matS^3$ or $\matB^3$;
\item the double projective plane $\twoRPtwo$, which is a Dehn surface of $\matRP^3$;
\item the four-hat, which is a Dehn surface of $\Lfourone$.
\end{itemize}
We conclude the proof by noting that in the first case $\Sigma'$ is derived from the sphere $\matS^2$.
\end{dimo}

\section{Subadditivity}\label{sec:subadditivity}

An important feature of a complexity function is to behave well with respect to the
cut-and-paste operations.
In this section, we will prove that the surface-complexity is subadditive under (boundary) connected sum.

\begin{teo}\label{teo:sub_additivity}
The surface-complexity of the connected sum and of the boundary connected sum of (connected and compact) 3-manifolds is less than or equal to the sum of their surface-complexities.
\end{teo}

\begin{dimo}
In order to prove the theorem, it is enough to prove the statement in the case where the number of the manifolds involved in the (boundary) connected sum is two.
Hence, if we call $M_1$ and $M_2$ the two manifolds, we need to prove that $sc(M_1\# M_2) \leqslant sc(M_1) + sc(M_2)$ and that $sc(M_1\#_{\partial} M_2) \leqslant sc(M_1) + sc(M_2)$.
Let $\Sigma_1$ (resp.~$\Sigma_2$) be a quasi-filling Dehn surface of $M_1$ (resp.~$M_2$) with $sc(M_1)$ (resp.~$sc(M_2)$) triple points, and let $\dot{M_1}$ (resp.~$\dot{M_2}$) be the punctured $M_1$ (resp.~$M_2$).

We start from the boundary connected sum.
It is obtained by identifying two discs $D_1$ and $D_2$ in $\partial M_1$ and $\partial M_2$, respectively.
Call $D$ the corresponding disc in $M_1\#_{\partial} M_2$, which is properly embedded.
Since the discs $D_1$ and $D_2$ are contained also in $\partial\dot{M_1}$ and $\partial\dot{M_2}$, respectively, we can consider also the (restriction of the) boundary connected sum
$\dot{M_1}\#_{\partial} \dot{M_2}$, which turns out to be a punctured $M_1\#_{\partial} M_2$.
We can suppose that $\Sigma_1$ and $\Sigma_2$ are embedded also in $\dot{M_1}\#_{\partial} \dot{M_2}$.
Consider an embedded arc $\alpha$ connecting $\Sigma_1$ and $\Sigma_2$ in $\dot{M_1}\#_{\partial} \dot{M_2}$ as shown in Fig.~\ref{fig:connected_sum}-left.
\begin{figure}
  \centerline{\includegraphics{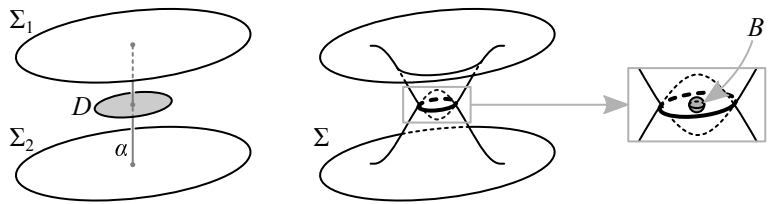}}
  \caption{The Dehn surface $\Sigma_1 \cup \Sigma_2$ in $\dot{M_1}\#_{\partial} \dot{M_2}$ with the arc $\alpha$ (left), its modification $\Sigma$ (centre), and the ball $B$ to be removed to get another punctured $M_1\#_{\partial} M_2$ (right).}
  \label{fig:connected_sum}
\end{figure}
The boundary connected sum $\dot{M_1}\#_{\partial} \dot{M_2}$ collapses to $\Sigma_1 \cup \Sigma_2 \cup \alpha$.
The last object is not a Dehn sphere (because it has a 1-dimensional part), but we can modify $\Sigma_1 \cup \Sigma_2$ as shown in Fig.~\ref{fig:connected_sum}-centre, getting a Dehn surface, say $\Sigma$.
If we remove a small ball $B$ (see Fig.~\ref{fig:connected_sum}-right) from $\dot{M_1}\#_{\partial} \dot{M_2}$ we obtain another punctured $M_1\#_{\partial} M_2$ which collapses to $\Sigma$.
Therefore, $\Sigma$ is a quasi-filling Dehn surface of $M_1\#_{\partial} M_2$ with $sc(M_1)+sc(M_2)$ triple points and hence we have $sc(M_1\#_{\partial} M_2) \leqslant sc(M_1) + sc(M_2)$.

For the connected sum we note that if we remove balls $B_1$ and $B_2$ from $M_1$ and $M_2$, respectively, and we apply a boundary connected sum along discs contained in $\partial B_1$ and $\partial B_2$, we get $M_1\# M_2$ minus a ball (whose boundary is the union of the complement of the discs in $\partial B_1$ and $\partial B_2$).
For $i=1,2$, up to applying a bubble-move (which does not increase the number of triple points), we can suppose that $\dot{M_i}$ is obtained from $M_i$ by removing at least one ball, and we can use one of the balls in $M_i\setminus\dot{M_i}$ as the ball $B_i$ for the construction of the connected sum.
By applying the construction above to the boundary connected sum along discs contained in $\partial B_1$ and $\partial B_2$, we get $\dot{M_1}\#_{\partial} \dot{M_2}$ and a Dehn surface $\Sigma$ to which $(\dot{M_1}\#_{\partial} \dot{M_2})\setminus B$ collapses (where $B$ is the small ball of Fig.~\ref{fig:connected_sum}-right).
Since $(\dot{M_1}\#_{\partial} \dot{M_2})\setminus B$ is a punctured $(M_1\setminus B_1)\#_{\partial} (M_2\setminus B_2)$ and the latter is $M_1\# M_2$ minus one ball, $(\dot{M_1}\#_{\partial} \dot{M_2})\setminus B$ is also a punctured $M_1\# M_2$.
Therefore, $\Sigma$ is a quasi-filling Dehn surface of $M_1\# M_2$ with $sc(M_1)+sc(M_2)$ triple points and hence we have $sc(M_1\# M_2) \leqslant sc(M_1) + sc(M_2)$.
\end{dimo}

\section{Estimations}

In general, calculating the surface-complexity $sc(M)$ of $M$ is very difficult, however it is relatively easy to estimate it.
More precisely, it is quite easy to give upper bounds for it.
Obviously, if we construct a quasi-filling Dehn surface $\Sigma$ of $M$, the number of triple points of $\Sigma$ is an upper bound for the surface-complexity of $M$.

The construction described in the first part of Theorem~\ref{teo:filling_surf_present}, obviously implies the following result.
\begin{teo}
If a (connected and compact) 3-manifold $M$ has an ideal triangulation with $n$ tetrahedra, the inequality $sc(M) \leqslant 4n$ holds.
\end{teo}

\begin{figure}
  \centerline{\includegraphics{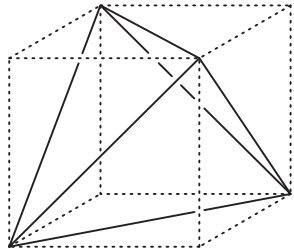}}
  \caption{Construction of an ideal triangulation from an ideal cubulation.}
  \label{fig:cube_to_tetra}
\end{figure}

\begin{figure}
  \centerline{\includegraphics{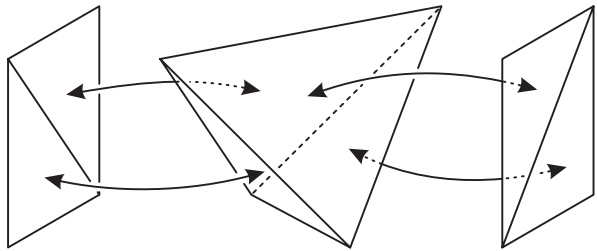}}
  \caption{Inserting a tetrahedron between two pairs of triangles not matching each other.}
  \label{fig:insert_tetra}
\end{figure}

\paragraph{From ideal cubulations to ideal triangulations}
We will now describe the inverse construction, allowing us to create ideal triangulations from filling Dehn surfaces (or ideal cubulations, by duality).
Let $\calC$ be an ideal cubulation of $M$.
Consider, for each cube of $\calC$, the five tetrahedra shown in Fig.~\ref{fig:cube_to_tetra}.
The idea is to glue together these ``bricks'' (each of which is made up of 5 tetrahedra) by following the identifications of the faces of $\calC$.
Note that each face of the cubes is divided by a diagonal into two triangles and that it may occur that these pairs of triangles do not match each other.
If this occurs, we insert a tetrahedron between them as shown in Fig.~\ref{fig:insert_tetra}.
Eventually, we get an ideal triangulation $\calT$ of $M$ with five tetrahedra for each cube of $\calC$ and at most one tetrahedron for each face of $\calC$.
Since the number of faces of an ideal cubulation is three times the number of cubes, we have that the number of tetrahedra in $\calT$ is at most 8 times the number of cubes.

We note that there are two different identifications of the abstract ``brick'' with each cube, so if there are $k$ cubes there are $2^k$ possibilities for the identifications with the cubes of $\calC$.
Some of them may need less insertions of tetrahedra (for matching the pairs of triangles in the faces of $\calC$) than others.
Hence, optimal choices may lead to a triangulation of $M$ with a number of tetrahedra closer to 5 times the number of cubes.

\paragraph{Matveev complexity}

The Matveev complexity~\cite{Matveev:compl} of a (connected and compact) 3-manifold $M$ is defined using simple spines.
A polyhedron $P$ is {\em simple} if the link of each point of $P$ can be embedded in the 1-skeleton of the tetrahedron.
The points of $P$ whose link is the whole 1-skeleton of the tetrahedron are called {\em vertices}.
If $\partial M\neq\emptyset$, a sub-polyhedron $P$ of $M$ is a {\em spine} of $M$ if $M$ collapses to $P$.
If $\partial M=\emptyset$, a sub-polyhedron $P$ of $M$ is a {\em spine} of $M$ if $M$ minus a ball collapses to $P$.
The {\em Matveev complexity} $c(M)$ of $M$ is the minimal number of vertices of a simple spine of $M$.
If $M$ is \ptwoirred, boundary-irreducible, without essential annuli and M\"obius strips, and different from $\matS^3$, $\matB^3$, $\matRP^3$ and $\Lthreeone$, its Matveev complexity is the minimal number of tetrahedra in an ideal triangulation of its (the Matveev complexity of $\matS^3$, $\matB^3$, $\matRP^3$ and $\Lthreeone$ is $0$), see~\cite{Matveev:book,Martelli-Petronio:decomposition}.
Therefore, the Matveev complexity is related to the surface-complexity, indeed the constructions described in Theorem~\ref{teo:filling_surf_present} and above, and the list of the \ptwoirred\ and boundary-irreducible 3-manifolds without essential annuli and M\"obius strips with Matveev complexy zero or surface-complexity zero obviously imply the following.

\begin{teo}
Let $M$ be a (connected and compact) \ptwoirred\ and boundary-irreducible 3-manifold without essential annuli and M\"obius strips, different from the lens spaces $\Lthreeone$ and $\Lfourone$; then the inequalities
$$
sc(M) \leqslant 4c(M)
\quad \mbox{and} \quad
c(M) \leqslant 8sc(M)
$$
hold.
Moreover, we have $c(\Lthreeone)=0$, $sc(\Lthreeone)>0$, $c(\Lfourone)>0$ and $sc(\Lfourone)=0$.
\end{teo}

\appendix

\section{The bidimensional case}

\begin{figure}
  \centerline{\includegraphics{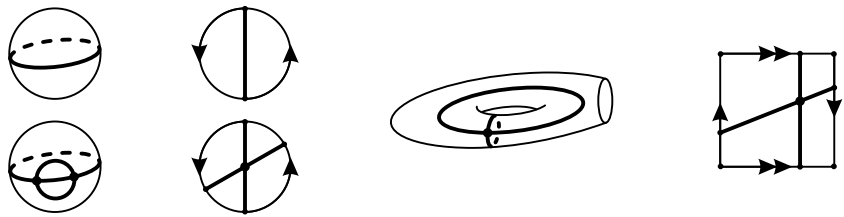}}
  \caption{Some Dehn loops of the sphere $\matS^2$, the projective plane $\matRP^2$, the torus minus a disc $\matT\setminus D$ and the Klein bottle $\matK$.}
  \label{fig:bidim_example}
\end{figure}

\begin{figure}
  \centerline{\includegraphics{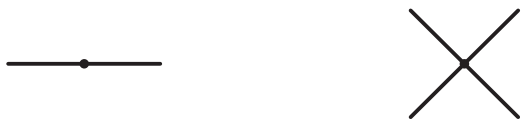}}
  \caption{Smooth points (left) and crossing points (right).}
  \label{fig:bidim_neigh}
\end{figure}

Let $S$ be a connected and compact surface.
Let us denote by $(\matS^1)^{\sqcup n}$ the disjoint union of $n$ circles $\matS^1$.
A subset $\Gamma$ of $S$ is said to be a {\em Dehn loop of $S$}
if there exists $n\in\matN$ and a transverse immersion $f\co (\matS^1)^{\sqcup n}\to S$ such that $\Gamma = f((\matS^1)^{\sqcup n})$.
Some examples are shown in Fig.~\ref{fig:bidim_example}.
There are only two types of points in $\Gamma$: smooth points and crossing points (see Fig.~\ref{fig:bidim_neigh}). The set of crossing points is denoted by $C(\Gamma)$.
As for 3-manifolds, we will need to remove some open discs from the surface $S$ such that their closures and the components of the boundary of $S$ are pairwise disjoint; we will denote the surface obtained by $\dot S$ and call it {\em a punctured $S$}.
A Dehn loop $\Gamma$ of $S$ will be called {\em quasi-filling} if a punctured $S$, $\dot S$, collapses to $\Gamma$.
A quasi-filling $\Gamma$ will be called {\em filling} if its singularities induce a cell-decomposition of $\Gamma$: more precisely,
\begin{itemize}
\item $C(\Gamma) \neq \emptyset$,
\item $\Gamma \setminus C(\Gamma)$ is made up of intervals.
\end{itemize}
Note that the latter condition is automatically satisfied once the former is fulfilled (because we are taking into account only connected surfaces).

Note that the Euler characteristic of $\dot S$ is equal to the Euler characteristic of the Dehn loop $\Gamma$.
On the other hand, if we start from an abstract Dehn loop $\Gamma$ and we thicken it to produce a surface with boundary, we can obtain different surfaces with boundary with the same Euler characteristic of $\Gamma$ by choosing the ``twists'' along the loops of $\Gamma$.
A simple case-by-case check implies that each surface $S$ has a filling Dehn loop.
Moreover, only the sphere $\matS^2$, the ball $\matB^2$, the annulus $\matA$, the projective plane $\matRP^2$ and the M\"obius strip $\matM$ have quasi-filling Dehn loops that are not filling, while all other surfaces have only filling Dehn loops.
However, as opposed to the 3-dimensional case (Theorem~\ref{teo:filling_surf_present}), a filling Dehn loop does not determine $S$ (even knowing the number of boundary components of $\dot S$ that must be filled up); for instance, the bouquet of two circles is a Dehn loop of the projective plane $\matRP^2$, the projective plane $\matRP^2$ with up to two discs removed, the torus $\matT$, the torus $\matT$ with one disc removed, the Klein bottle $\matK$ and the Klein bottle $\matK$ with one disc removed (see Fig.~\ref{fig:bidim_example}).

We define the {\em loop-complexity} $lc(S)$ of a connected and compact surface $S$ as the minimal number of crossing points of a quasi-filling Dehn loop of $S$.
For a non-closed $S$ different from $\matB^2$, $\matA$ and $\matM$, by means of an easy Euler characteristic argument (in dimension 2 it is a very powerful tool) and a case-by-case analysis, it is quite easy to prove that the minimal number of crossing points of a filling Dehn surface of $S$ occurs when there is no puncture (i.e.~$S=\dot S$) and that this number of crossing points is the opposite of the Euler characteristic of $S$.
Since for closed surfaces a puncture is needed (to create the boundary for collapsing), we obtain that the loop-complexity of a surface $S$ with Euler characteristic $\chi$ is
\begin{itemize}
\item
$-\chi$ if $\partial S\neq\emptyset$, except for $\matB^2$ having loop-complexity 0;
\item
$1-\chi$ if $\partial S=\emptyset$, except for $\matS^2$ having loop-complexity 0.
\end{itemize}
Note that it is not true that the loop-complexity of the connected sum of two surfaces is at most the sum of their loop-complexities; for instance, we have $K=\matRP^2\#\matRP^2$ while $lc(K)=1 \not\leqslant 0=lc(\matRP^2)+lc(\matRP^2)$.

We can also consider, as done for 3-manifolds, {\em ideal cubulations} $\calC$ of connected and compact surfaces $S$, i.e.~cell-decompositions of the topological space $\widehat S$ obtained from $\dot S$ by collapsing each boundary component (of $\dot S$) to a point such that each 2-cell (a {\em square}) is glued along four edges.
An ideal cubulation of $S$ can be constructed from a filling Dehn loop $\Gamma$ of $S$ by considering an abstract square for each crossing point of $\Gamma$ and by gluing the squares together along the edges.
However, there are two possibilities for gluing two squares along an edge and the abstract polyhedron $\Gamma$ does not encode any information to choose the right one.
(In some sense, this explains why a filling Dehn loop does not determine unambiguously one surface.)
If we consider also the immersion of $(\matS^1)^{\sqcup n}$ in the surface $S$ containing the Dehn loop $\Gamma$, we can choose the right identifications and construct an ideal cubulation $\calC$ of $S$.

The converse construction can also be performed, obtaining a filling Dehn loop of a connected and compact surface $S$ from an ideal cubulation of $S$.
These constructions tell us that $lc(S)$ is the minimal number of squares in a cubulation of $S$, except for $\matS^2$, $\matB^2$, $\matA$, $\matRP^2$ and $\matM$, whose loop-complexity is $0$.

\begin{small}

\end{small}


\begin{thebibliography}{99}

\bibitem{Aitchison-Matsumotoi-Rubinstein}
\textsc{I.~R.~Aitchison -- S.~Matsumotoi -- J.~H.~Rubinstein},
\textit{Immersed surfaces in cubed manifolds},
Asian J. Math. {\bf 1} (1997), no.~1, 85--95.

\bibitem{Amendola:surf_inv}
\textsc{G.~Amendola},
\textit{A Local Calculus for Nullhomotopic Filling Dehn Spheres},
Algebr. Geom. Topol. {\bf 9} (2009), no.~2, 903--933.

\bibitem{Amendola:surf_compl}
\textsc{G.~Amendola},
\textit{A 3-Manifold Complexity via Immersed Surfaces},
J. Knot Theory Ramif. {\bf 19} (2010), no.~12, 1549--1569.

\bibitem{Amendola:sc1}
\textsc{G.~Amendola},
\textit{Orientable closed 3-manifolds with surface-complexity one},
Atti Semin. Mat. Fis. Univ. Modena Reggio Emilia {\bf 57} (2010), 17--26.

\bibitem{Babson-Chan}
\textsc{E.~K.~Babson -- C.~S.~Chan},
\textit{Counting faces of cubical spheres modulo two},
Discrete Math. {\bf 212} (2000), no.~3, 169--183.

\bibitem{Casler}
\textsc{B.~G.~Casler},
\textit{An imbedding theorem for connected $3$-manifolds with boundary},
Proc. Amer. Math. Soc. {\bf 16} (1965), 559--566.

\bibitem{Dolbilin-Shtanko-Shtogrin}
\textsc{N.~P.~Dolbilin -- M.~A.~Shtan$'$ko -- M.~I.~Shtogrin},
\textit{Cubic manifolds in lattices},
Izv. Ross. Akad. Nauk Ser. Mat. {\bf 58} (1994), no.~2, 93--107;
translation in 
Russian Acad. Sci. Izv. Math. {\bf 44} (1995), no.~2, 301--313.

\bibitem{Fenn-Rourke}
\textsc{R.~Fenn -- C.~Rourke},
\textit{Nice spines of 3-manifolds},
Topology of low-dimensional manifolds, Lecture Notes in Math., vol. 722 (1977), 31–36, Springer, Berlin 1979.

\bibitem{Funar}
\textsc{L.~Funar},
\textit{Cubulations, immersions, mappability and a problem of Habegger},
Ann. Sci. \'Ecole Norm. Sup. (4) {\bf 32} (1999), no.~5, 681--700.

\bibitem{Korablev-Kazakov:cubic_complexity_two}
\textsc{F.~G.~Korablev -- A.~A.~Kazakov},
\textit{Manifolds of cubic complexity two},
Sib. \`Elektron. Mat. Izv. {\bf 13} (2016), 1–15.

\bibitem{Lozano-Vigara:subadditivity}
\textsc{A.~Lozano -- R.~Vigara},
\textit{On the subadditivity of Montesinos complexity of closed orientable 3-manifolds},
Rev. R. Acad. Cienc. Exactas Fís. Nat. Ser. A Mat. RACSAM {\bf 109} (2015), no.~2, 267–279. 

\bibitem{Lozano-Vigara:representing}
\textsc{A.~Lozano -- R.~Vigara},
``Representing 3-manifolds by filling Dehn surfaces,''
Series on Knots and Everything, 58. World Scientific Publishing Co. Pte. Ltd., Hackensack, NJ, 2016.
xvii+276 pp.

\bibitem{Lozano-Vigara:triple_point_spectrum}
\textsc{A.~Lozano -- R.~Vigara},
\textit{The Triple-Point Spectrum of Closed Orientable 3-Manifolds},
Mediterranean Journal of Mathematics {\bf 16} (2019), no.~2, Paper No. 71, 19 pp.

\bibitem{Martelli-Petronio:decomposition}
\textsc{B.~Martelli -- C.~Petronio},
\textit{A new decomposition theorem for 3-manifolds},
Ill. J. Math. {\bf 46} (2002), no.~3, 755–780

\bibitem{Matveev:compl}
\textsc{S.~V.~Matveev},
\textit{The theory of the complexity of three-dimensional manifolds},
Akad. Nauk Ukrain. SSR Inst. Mat. Preprint (1988), no.~13, 32~pp.

\bibitem{Matveev:book}
\textsc{S.~V.~Matveev},
``Algorithmic topology and classification of 3-manifolds,''
Algorithms and Computation in Mathematics, 9. Springer-Verlag, Berlin, 2003.
xii+478 pp.

\bibitem{Montesinos}
\textsc{J.~M.~Montesinos-Amilibia},
\textit{Representing 3-manifolds by Dehn spheres},
Mathematical contributions: volume in honor of Professor Joaqu\'\i n Arregui Fern\'andez, 239--247, Homen. Univ. Complut., Editorial Complutense, Madrid, 2000.

\bibitem{Papa}
\textsc{C.~D.~Papakyriakopoulos},
\textit{On Dehn's lemma and the asphericity of knots},
Ann. of Math. (2) {\bf 66} (1957), 1--26.

\bibitem{Rourke}
\textsc{C.~Rourke},
\textit{A new proof that $\Omega_3$ is zero},
J. London Math. Soc. (2) {\bf 31} (1985), no. 2, 373--376.

\bibitem{Shtanko-Shtogrin}
\textsc{M.~A.~Shtan$'$ko -- M.~I.~Shtogrin},
\textit{Embedding cubic manifolds and complexes into a cubic lattice},
Uspekhi Mat. Nauk {\bf 47} (1992), no.~1 (283), 219--220;
translation in 
Russian Math. Surveys {\bf 47} (1992), no.~1, 267--268.

\bibitem{Vigara:present}
\textsc{R.~Vigara},
\textit{A new proof of a theorem of J.~M.~Montesinos},
J. Math. Sci. Univ. Tokyo {\bf 11} (2004), no. 3, 325--351.

\bibitem{Vigara:calculus}
\textsc{R.~Vigara},
\textit{A set of moves for Johansson representation of 3-manifolds},
Fund. Math. {\bf 190} (2006), 245--288.

\bibitem{Vigara:lifting}
\textsc{R.~Vigara},
\textit{Lifting filling Dehn spheres},
J. Knot Theory Ramifications {\bf 21} (2012), no. 8, 1250082, 7 pp.

\end{thebibliography}
\end{document}